\documentclass{elsarticle}

\usepackage{amsmath,amssymb,amsthm}
\usepackage{graphicx}

\usepackage{lineno,hyperref}
\modulolinenumbers[5]

\journal{Applied Mathematics and Computation}

\newtheorem{theorem}{Theorem}

\begin{document}

\begin{frontmatter}

\title{A new approach to predict changes in physical condition: A new
extension of the classical Banister model}

\author{Marcos Matabuena,  Rosana Rodr\'iguez-L\'opez}
\address{Departamento de An\'alise Matem\'atica, Estat\'istica e Optimizaci\'on,\\
 Universidade de Santiago de Compostela,\\
15782, Santiago de Compostela, Spain}
\ead{matabuelas@yahoo.es, rosana.rodriguez.lopez@usc.es}

\begin{abstract}
In this article, a new model based on techniques of differential equations is introduced to predict the athletic performance based training load and a data sample of the physical
form of athletes arises. This model is an extension of the classical model of Banister but, in this case, unlike the classical Banister model, the
variation produced in the athletic performance depends, not only on the current training load,
but also on the training performed the previous day. 
The model
has been validated with the training data of a cyclist taken from the reference 
\cite{Clarke}, obtaining an excellent fit of the predicted data with respect to the experimental data.
\end{abstract}


\end{frontmatter}

\linenumbers

% main text

\section{Introduction}
Athletes, in their quest to optimize performance, try to use mathematical tools
(Issurin, 2010 \cite{Issurin}), but nowadays there are no tools to know the effect of workout
on the final performance (Gonz\'alez Badillo \& Gorostiaga, 1995 \cite{GonzBad}).  A first attempt
at the problem of modeling and predicting changes induced by fitness training
appeared in 1975 with the contribution of the physiologist Banister and his coworkers (Banister
et al., 1975 \cite{Banister75}), whose model has been frequently used in the literature to predict
the effects of training on physical performance, see, for example, its validation
in triathlon athletes by Banister et al. in 1999 \cite{Banister99}, the case of a swimmer presented in 1976 by Calvert et al.  \cite{Calvert}, studies in cross-country skiers made in 1992 by Candau et al.  \cite{Candau}, the convenience to
exercise physiology practitioners and researchers in the use of mathematical models
(see Clarke \& Skiba, 2013 \cite{Clarke}), the optimal design of a training strategy (see Fitz-Clarke
et al., 1991 \cite{Fitz}), the performance of runners analyzed in 1990 by Morton et al.  \cite{Morton1}, or the recent consideration of altitude training by Rodr\'iguez et al. \cite{Rodr}, among others.
The interest of the application of computational intelligence in sports has also been defended
(see \cite{AMC}).

Since the introduction of the Banister model, there have been proposed
modified versions of the model, some of which are those introduced by Busso in
2003 \cite{Busso}, or Coggan in 2006 \cite{Coggan}, among others, and some new
models have been built through regression techniques, neural networks (see Pfeiffer
\& Hohmann, 2012 \cite{PH}) or time series (see Pfeiffer, 2008).
In this paper, we propose a new model based on techniques of functional
differential equations and nonlinear regression that substantially improves the
accuracy of the above methods, besides overcoming some of their limitations
(see Hellard et al., 2006 \cite{Hellard}). The model proposed has been successfully applied to
study and predict the effect of training a cyclist in the task of optimizing a
workout plan, constituting what might be a reliable and useful tool for the future, in
their decision making and when applied to other sports.

\section{Mathematical model}
In the work \cite{art},  the authors consider the following equation:
$$g'(t)=-\frac{1}{\tau_1}g(t)+w(t),$$
whose solution is given by
$$g(t)=g(0)+\int_{0}^t w(s) e^{\frac{-1}{\tau_1}(t-s)}\, ds=\int_{0}^t w(s) e^{\frac{-1}{\tau_1}(t-s)}\, ds, $$
which gives the following approximation
$$g(n)=\sum_{i=0}^{n-1} w(i) e^{\frac{-1}{\tau_1}(n-i)} .$$
If we consider  $w(0)=0$, then it coincides with  equation (3) in
\cite{art}.

\bigskip

Next, we consider the delay differential model
\begin{equation}
\label{1} g'(t)=-\frac{1}{\tau_1}g(t)-\frac{1}{\tau_2} g(t-1)+w(t).
\end{equation}
and consider an equivalent formulation.
\begin{theorem}\label{teo1}
The equation (\ref{1})
 is equivalent to 
$$g(t) e^{\frac{1}{\tau_1}(t-k)}-g(k)=\int_{k}^t \left[w(s) -\frac{1}{\tau_2} g(s-1)\right] e^{\frac{1}{\tau_1}(s-k)}\, ds,$$
for $t\in (k,k+1]$, $ k=0,1,\ldots$
\end{theorem}

From Theorem \ref{teo1}, we have
$$g(t) =g(k) e^{-\frac{1}{\tau_1}(t-k)}+ \int_{k}^t \left[w(s) -\frac{1}{\tau_2} g(s-1)\right] e^{-\frac{1}{\tau_1}(t-s)}\, ds,$$
for $t\in (k,k+1]$, $k=0,1,\ldots $,
so that we can consider the following approximation:
\begin{equation*}
\begin{split}
g(k+1) &=g(k) e^{-\frac{1}{\tau_1}}+  \left[w(k) -\frac{1}{\tau_2} g(k-1)\right] e^{-\frac{1}{\tau_1}}\\
&=
\left[w(k)+g(k) -\frac{1}{\tau_2} g(k-1)\right] e^{-\frac{1}{\tau_1}}
, k=0,1,\ldots 
\end{split}
\end{equation*}
Note that, for $t \in [0,1]$, we get
$$g(t) = \int_{0}^t \left[w(s) -\frac{1}{\tau_2} g(s-1)\right] e^{-\frac{1}{\tau_1}(t-s)}\, ds= \int_{0}^t w(s)  e^{-\frac{1}{\tau_1}(t-s)}\, ds,
$$
which is coincident with \cite{art}.

\bigskip

Other way to calculate the solution to \eqref{1}, by integration between  $0$ and $t$, which
provides a model  more similar to that in \cite{art},  would be to write the equation
 equivalently as
$$g(t) e^{\frac{1}{\tau_1}t}=g(0)+\int_{0}^t \left[w(s) -\frac{1}{\tau_2} g(s-1)\right] e^{\frac{1}{\tau_1}s}\, ds, \quad  t\geq 0,$$
or
\begin{equation*}
\begin{split}
g(t)&=g(0) e^{-\frac{1}{\tau_1}t}+\int_{0}^t \left[w(s) -\frac{1}{\tau_2} g(s-1)\right] e^{-\frac{1}{\tau_1}(t-s)}\, ds\\
& =\int_{0}^t \left[w(s) -\frac{1}{\tau_2} g(s-1)\right] e^{-\frac{1}{\tau_1}(t-s)}\, ds, \quad  t\geq 0,
\end{split}
\end{equation*}
therefore 
$$g(n)=\int_{0}^n \left[w(s) -\frac{1}{\tau_2} g(s-1)\right] e^{-\frac{1}{\tau_1}(n-s)}\, ds, \quad  t\geq 0,$$
and we can take, similarly to the study in \cite{art}, the approximation
$$g(n)=\sum_{i=0}^{n-1} \left[w(i)-\frac{1}{\tau_2} g(i-1)\right]  e^{\frac{-1}{\tau_1}(n-i)} .$$
If we suppose that  $w(0)=0$, since $g(-1)=0$, the previous expression is equal to
\begin{equation}
\label{difer}
g(n)=\sum_{i=1}^{n-1} \left[w(i)-\frac{1}{\tau_2} g(i-1)\right]  e^{\frac{-1}{\tau_1}(n-i)} .
\end{equation}

\subsection{Model proposed}

The model proposed is based on the study of the properties of differential equations with delay and the estimation of parameters. We first state the problem in mathematical terms, as follows:
\begin{equation}\label{ec1}
p(t)=p(0)+k_1g(t)-k_2h(t),
\end{equation}
where $k_1$ and   $k_2$ are two positive constants and  $p(0)$ is the initial physical condition. The values of the functions $g(t)$ and $h(t)$  are obtained by solving the following system of differential equations:
$$
\left\{\begin{array}{c}
g'(t)=w(t)-\frac{1}{\tau_1}g(t)-\frac{1}{\tau_2}g(t-1),\\ \\
h'(t)=w(t)-\frac{1}{\tau_3}h(t)-\frac{1}{\tau_4}h(t-1),
\end{array}\right.$$
where $w(t)$ represents the load at the instant of time $t$ and $\tau_1$, $\tau_2$, $\tau_3$, $\tau_4$  are positive constants. Since, in practice, we work with measurements at fixed times, we consider a discrete model and approximate the solution to the differential system by considering the following difference equation (see \eqref{difer}):
$$p(n)=k_1 \sum_{i=1}^{n-1} \left[  w(i)-\frac{1}{\tau_2}g(i-1)\right]e^{\frac{1}{\tau_1}(n-i)}
-k_2 \sum_{i=1}^{n-1} \left[  w(i)-\frac{1}{\tau_4}h(i-1)\right]e^{\frac{1}{\tau_3}(n-i)},
$$
where $n$ denotes the $n$th instant of time. The last equation represents the explicit solution of the  discretized (difference) equation of problem \eqref{ec1}.

At this point, we can use nonlinear regression techniques, or any heuristic algorithm, to estimate the constants  $k_1, k_2, \tau_1,\tau_2,\tau_3,\tau_4$ by the method of least squares, using a sample of physical condition $p(1),\ldots, p(s)$, at different times. In this work, we use the experimental data of a cyclist extracted from the 2013 article by Clarke \& Skiba \cite{Clarke}. To optimize the constants, we use the Solver in Microsoft Office and we represent the results with the statistical package R.

\subsection{Results obtained}
The result of the $R^2$ test between the experimental data and the mathematical model proposed is $0.9903$. In Figure 1, we find a graph where the red points represent the data and the blue curve the solution predicted by the model. In Figure 2, we see the athlete on his training-workload, which determines the change in his physical condition. We observe that this graph presents high and low values, which correspond to the effect of the assimilation
of the training and the fatigue states, respectively, what actually happens in reality.

\begin{figure}
\centering
\includegraphics[width=9cm]{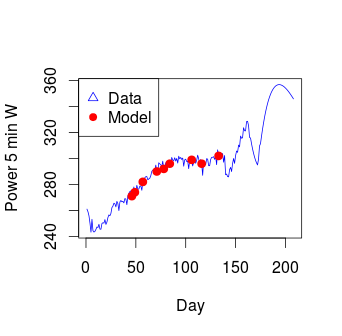}
\caption{Experimental data (in red) and the solution predicted by the model (in blue)}
\end{figure}

\begin{figure}
	\centering
	\includegraphics[width=8cm]{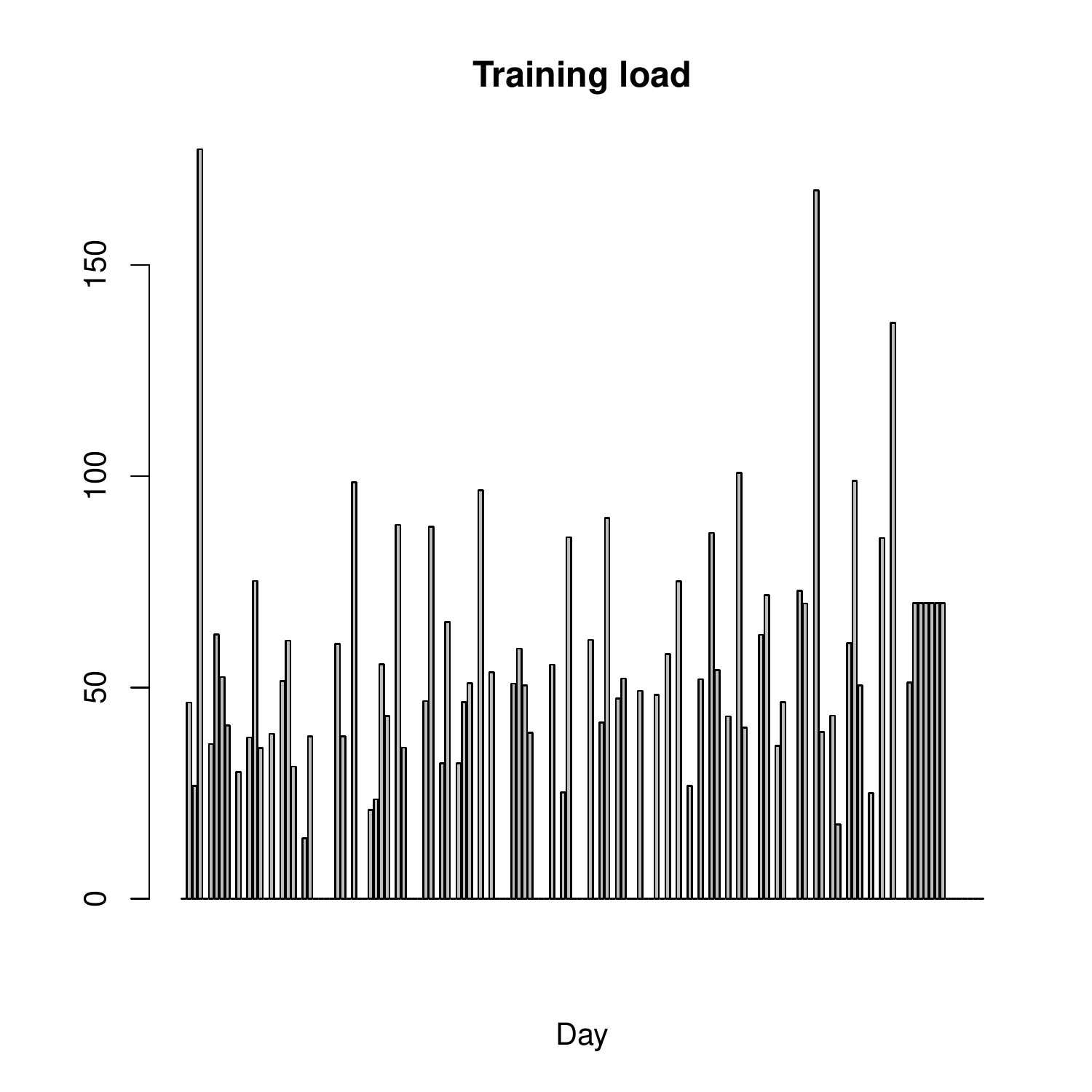}
	\caption{Function $w(t)$}
\end{figure}

\section{Other formulations of the model}
In this section, we show other different formulations of the model.

\subsection{Problem with three delays}
Next, we consider
\begin{equation}
\label{2} g'(t)=-\frac{1}{\tau_1}g(t)-\frac{1}{\tau_2} g(t-1)-\frac{1}{\tau_3} g(t-2)-\frac{1}{\tau_4} g(t-3)+w(t).
\end{equation}
The solution to (\ref{2}) can be written as:
$$g(t)=g(0) e^{-\frac{1}{\tau_1}t}+\int_{0}^t \left[w(s) -\sum_{k=1}^3 \frac{1}{\tau_{k+1}} g(s-k)\right] e^{-\frac{1}{\tau_1}(t-s)}\, ds, \quad  t\geq 0,$$
that is,
$$
g(t)=\int_{0}^t \left[w(s) -\sum_{k=1}^3 \frac{1}{\tau_{k+1}} g(s-k)\right]  e^{-\frac{1}{\tau_1}(t-s)}\, ds, \quad  t\geq 0,$$
so that, 
$$g(n)=\int_{0}^n 
\left[w(s) -\sum_{k=1}^3 \frac{1}{\tau_{k+1}} g(s-k)\right]
 e^{-\frac{1}{\tau_1}(n-s)}\, ds$$
and we can take, similarly to \cite{art}, the approximation
$$g(n)=\sum_{i=0}^{n-1} \left [w(i) -\sum_{k=1}^3 \frac{1}{\tau_{k+1}} g(i-k)\right]  e^{\frac{-1}{\tau_1}(n-i)} $$
or
$$g(n)
=\sum_{i=0}^{n-1} \left [w(i) - \frac{1}{\tau_2} g(i-1)-\frac{1}{\tau_3} g(i-2)-\frac{1}{\tau_4} g(i-3)\right]  e^{\frac{-1}{\tau_1}(n-i)} 
.$$
If we suppose that $w(0)=0$, since $g(s)=0$, $\forall s\in [-3,0]$, we have
\begin{equation*}
\begin{split}
g(n) &
=\sum_{i=1}^{n-1} \left [w(i) - \frac{1}{\tau_2} g(i-1)-\frac{1}{\tau_3} g(i-2)-\frac{1}{\tau_4} g(i-3)\right]  e^{\frac{-1}{\tau_1}(n-i)}\\
& =\sum_{i=1}^{n-1} w(i) e^{\frac{-1}{\tau_1}(n-i)}- \frac{1}{\tau_2} \sum_{i=1}^{n-1} g(i-1)e^{\frac{-1}{\tau_1}(n-i)}
-\frac{1}{\tau_3} \sum_{i=1}^{n-1}g(i-2)e^{\frac{-1}{\tau_1}(n-i)}\\
& \hspace*{3em}
-\frac{1}{\tau_4} \sum_{i=1}^{n-1}g(i-3) e^{\frac{-1}{\tau_1}(n-i)}\\
& =\sum_{i=1}^{n-1} w(i) e^{\frac{-1}{\tau_1}(n-i)}- \frac{1}{\tau_2} \sum_{i=1}^{n-1} g(i-1)e^{\frac{-1}{\tau_1}(n-i)}
-\frac{1}{\tau_3} \sum_{j=0}^{n-2}g(j-1)e^{\frac{-1}{\tau_1}(n-j-1)}\\
&\hspace*{3em}
-\frac{1}{\tau_4} \sum_{j=-1}^{n-3}g(j-1) e^{\frac{-1}{\tau_1}(n-j-2)}\\
& =\sum_{i=1}^{n-1} w(i) e^{\frac{-1}{\tau_1}(n-i)}- \frac{1}{\tau_2} \sum_{i=1}^{n-1} g(i-1)e^{\frac{-1}{\tau_1}(n-i)}
-\frac{1}{\tau_3} \sum_{j=1}^{n-2}g(j-1)e^{\frac{-1}{\tau_1}(n-j-1)}\\
&\hspace*{3em}
-\frac{1}{\tau_4} \sum_{j=1}^{n-3}g(j-1) e^{\frac{-1}{\tau_1}(n-j-2)}\\
& =\sum_{i=1}^{n-1} w(i) e^{\frac{-1}{\tau_1}(n-i)}
-\sum_{i=1}^{n-3} g(i-1) \left[ \frac{1}{\tau_2} e^{\frac{-1}{\tau_1}(n-i)} 
+\frac{1}{\tau_3} e^{\frac{-1}{\tau_1}(n-i-1)}
+\frac{1}{\tau_4} e^{\frac{-1}{\tau_1}(n-i-2)}
 \right]\\
 &\hspace*{3em} -g(n-3) \left[  \frac{1}{\tau_2} e^{\frac{-2}{\tau_1}}
+\frac{1}{\tau_3} e^{\frac{-1}{\tau_1}}  \right]
-g(n-2) \frac{1}{\tau_2} e^{\frac{-1}{\tau_1}}\\
& =\sum_{i=1}^{n-1} w(i) e^{\frac{-1}{\tau_1}(n-i)}
-\sum_{i=1}^{n-3} g(i-1) \left[ \frac{1}{\tau_2}
+\frac{1}{\tau_3} e^{\frac{1}{\tau_1}}
+\frac{1}{\tau_4} e^{\frac{2}{\tau_1}}
 \right]  e^{\frac{-1}{\tau_1}(n-i)}  \\
 &\hspace*{3em} -g(n-3) \left[  \frac{1}{\tau_2} e^{\frac{-2}{\tau_1}}
+\frac{1}{\tau_3} e^{\frac{-1}{\tau_1}}  \right]
-g(n-2) \frac{1}{\tau_2} e^{\frac{-1}{\tau_1}}.
\end{split}
\end{equation*}

Following for problem (\ref{2}) a formulation similar to the first one exposed,  we have:
$$g(t) =g(k) e^{-\frac{1}{\tau_1}(t-k)}+ \int_{k}^t \left[w(s) -\frac{1}{\tau_2} g(s-1)-\frac{1}{\tau_3} g(s-2)
-\frac{1}{\tau_4} g(s-3)\right] e^{-\frac{1}{\tau_1}(t-s)}\, ds,$$
for $t\in (k,k+1], $ $k=0,1,\ldots $,
so that, we can give the following approximation:
\begin{equation}\label{recurr2}
\begin{split}
g(k+1) &=g(k) e^{-\frac{1}{\tau_1}}+  \left[w(k) -\frac{1}{\tau_2} g(k-1)-\frac{1}{\tau_3} g(k-2)
-\frac{1}{\tau_4} g(k-3)\right] e^{-\frac{1}{\tau_1}}\\
& =
\left[w(k)+g(k) -\frac{1}{\tau_2} g(k-1)-\frac{1}{\tau_3} g(k-2)
-\frac{1}{\tau_4} g(k-3)\right] e^{-\frac{1}{\tau_1}}
,
\end{split}
\end{equation}
for $ k=0,1,\ldots $

\subsection{Modelling through an integral equation}
We consider
\begin{equation}
\label{3} g'(t)=-\frac{1}{\tau_1}g(t)+\tau_5 \int_0^t \delta(t,s)g(s)\, ds+w(t).
\end{equation}
We can impose, for instance, that
$\delta(t,s)=0$ for $s<t-4$.
The solution to (\ref{3}) is obtained as
\begin{equation*}
\begin{split}
g(t)&=g(0) e^{-\frac{1}{\tau_1}t}+\int_{0}^t \left[w(u) +\tau_5 \int_0^u \delta(u,s)g(s)\, ds\right] e^{-\frac{1}{\tau_1}(t-u)}\, du\\
&=\int_{0}^t \left[w(u) +\tau_5 \int_0^u \delta(u,s)g(s)\, ds\right] e^{-\frac{1}{\tau_1}(t-u)}\, du, \quad  t\geq 0,
\end{split}
\end{equation*}
what gives the approximation
\begin{equation*}
\begin{split}
g(n)&=\sum_{i=0}^{n-1} \left[w(i) +\tau_5 \int_0^i \delta(i,s)g(s)\, ds\right] e^{-\frac{1}{\tau_1}(n-i)}\\
& =\sum_{i=0}^{n-1} \left[w(i) +\tau_5 \sum_{j=0}^{i-1} \delta(i,j)g(j)\right] e^{-\frac{1}{\tau_1}(n-i)}.
\end{split}
\end{equation*}
From this identity, taking into account that $w(0)=0$ and $g(0)=0$, we get
\begin{equation*}
\begin{split}
g(n)& =\sum_{i=0}^{n-1} w(i) e^{-\frac{1}{\tau_1}(n-i)}  +\tau_5\sum_{i=0}^{n-1}  \sum_{j=0}^{i-1} \delta(i,j)g(j) e^{-\frac{1}{\tau_1}(n-i)}\\
& =\sum_{i=1}^{n-1} w(i) e^{-\frac{1}{\tau_1}(n-i)}  +\tau_5\sum_{i=0}^{n-1}  \sum_{j=1}^{i-1} \delta(i,j)g(j) e^{-\frac{1}{\tau_1}(n-i)}\\
\end{split}
\end{equation*}
If we understand that $\delta(i,j)\neq 0$ only if  $j=i-3,i-2,i-1$, 
we get 
\begin{equation*}
\begin{split}
& \sum_{i=0}^{n-1}  \sum_{j=1}^{i-1} \delta(i,j)g(j) e^{-\frac{1}{\tau_1}(n-i)}\\&
=   
  \delta(2,1)g(1) e^{-\frac{1}{\tau_1}(n-2)}+  \sum_{j=1}^{2} \delta(3,j)g(j) e^{-\frac{1}{\tau_1}(n-3)}  +
   \sum_{j=1}^{3} \delta(4,j)g(j) e^{-\frac{1}{\tau_1}(n-4)}\\
   &\hspace*{2em}+
    \sum_{j=1}^{4} \delta(5,j)g(j) e^{-\frac{1}{\tau_1}(n-5)}+\cdots+ \sum_{j=1}^{n-3} \delta(n-2,j)g(j) e^{-\frac{2}{\tau_1}}\\
    & \hspace*{2em} + \sum_{j=1}^{n-2} \delta(n-1,j)g(j) e^{-\frac{1}{\tau_1}}\\
    & 
=   
  \delta(2,1)g(1) e^{-\frac{1}{\tau_1}(n-2)}+  \sum_{j=1}^{2} \delta(3,j)g(j) e^{-\frac{1}{\tau_1}(n-3)}  +
   \sum_{j=1}^{3} \delta(4,j)g(j) e^{-\frac{1}{\tau_1}(n-4)}\\
   &\hspace*{2em} +
    \sum_{j=2}^{4} \delta(5,j)g(j) e^{-\frac{1}{\tau_1}(n-5)}+\cdots+ \sum_{j=n-5}^{n-3} \delta(n-2,j)g(j) e^{-\frac{2}{\tau_1}}\\
    & \hspace*{2em} + \sum_{j=n-4}^{n-2} \delta(n-1,j)g(j) e^{-\frac{1}{\tau_1}}\\
    & = \left[ \delta(2,1) e^{-\frac{1}{\tau_1}(n-2)}+  \delta(3,1) e^{-\frac{1}{\tau_1}(n-3)}  
    + \delta(4,1) e^{-\frac{1}{\tau_1}(n-4)}\right] g(1)\\
    & \hspace*{2em} +
    \left[   \delta(3,2) e^{-\frac{1}{\tau_1}(n-3)}  +\delta(4,2) e^{-\frac{1}{\tau_1}(n-4)} +
     \delta(5,2) e^{-\frac{1}{\tau_1}(n-5)}   \right] g(2)+
\cdots \\
&\hspace*{2em} + \left[  \delta(n-3,n-4) e^{-\frac{3}{\tau_1}} +  \delta(n-2,n-4) e^{-\frac{2}{\tau_1}} 
+  \delta(n-1,n-4) e^{-\frac{1}{\tau_1}}   \right] g(n-4)\\
   &\hspace*{2em} +\left[  \delta(n-2,n-3) e^{-\frac{2}{\tau_1}} 
+  \delta(n-1,n-3) e^{-\frac{1}{\tau_1}}   \right] g(n-3)\\
& \hspace*{2em}
+  \delta(n-1,n-2) e^{-\frac{1}{\tau_1}}g(n-2)\\
&  =\sum_{j=1}^{n-4} \left[ \sum_{i=j+1}^{j+3} \delta(i,j) e^{-\frac{1}{\tau_1}(n-i)}  \right] g(j)
+ \left[ \sum_{i=n-2}^{n-1} \delta(i,n-3) e^{-\frac{1}{\tau_1}(n-i)}\right]   g(n-3)\\
& \hspace*{2em}
+  \delta(n-1,n-2) e^{-\frac{1}{\tau_1}}g(n-2).
\end{split}
\end{equation*}

We can proceed by assigning  fixed values in $(0,1)$ to
$$\delta(i,i-1),\:  \delta(i,i-2),\:  \delta(i,i-3),$$
that can be always of the same type, for instance,
$$\delta(i,i-1)=0.5,\:  \delta(i,i-2)=0.3,\:  \delta(i,i-3)=0.2.$$

Following for problem (\ref{3}) a formulation similar to the first one exposed, we get:
$$g(t) =g(k) e^{-\frac{1}{\tau_1}(t-k)}+ \int_{k}^t  \left[w(u) +\tau_5 \int_0^u \delta(u,s)g(s)\, ds\right] e^{-\frac{1}{\tau_1}(t-u)}\, du,$$
for $t\in (k,k+1], $ $k=0,1,\ldots $,
so that, we can give the approximation:
$$
g(k+1) =g(k) e^{-\frac{1}{\tau_1}}+  \left[w(k)  +\tau_5 \int_0^k \delta(k,s)g(s)\, ds \right] e^{-\frac{1}{\tau_1}}$$
and, taking again an approximation,  we get the recursive formula
$$
g(k+1) =
\left[w(k)+g(k) +\tau_5 \sum_{j=k-3}^{k-1} \delta(k,j)g(j) \right]  e^{-\frac{1}{\tau_1}}
, k=0,1,\ldots 
$$
If we select always
$$\delta(k,k-1)=0.5,\:  \delta(k,k-2)=0.3,\:  \delta(k,k-3)=0.2,$$
the recurrence formula is given by:
\begin{equation}\label{recurr3}
g(k+1) =
\left[w(k)+g(k) +0.5\,  \tau_5 \, g(k-1)+ 0.3\, \tau_5 \, g(k-2)+ 0.2 \,\tau_5 \, g(k-3)  \right]  e^{-\frac{1}{\tau_1}}
, 
\end{equation}
for $k=0,1,\ldots$, 
so that, choosing  $\delta$ this way, we obtain something similar  to the retarded approach,
since choosing  $\tau_2$, $\tau_3$ and $\tau_4$ such that: $-\frac{1}{\tau_2}=0.5 \tau_5$, 
$-\frac{1}{\tau_3}=0.3 \tau_5$ and $-\frac{1}{\tau_4}=0.2 \tau_5$, the recurrence formula
(\ref{recurr3}) coincides with (\ref{recurr2}).

\section{Conclusions}
In this work, we have shown that, using the model proposed, it is possible to obtain an excellent fit with respect to the experimental data, calculating the positive and negative effects of the training with the ability to predict the effect of training loads in fitness. The improvement of the accuracy in the predicted data leads us definitely to an important tool for planning sports coaching, and new gadgets and software that attempt to control the effect that sports training has on the athlete individually. Besides, this model could be adapted for weight control or to predict the risk of injury or the emotional state of the athlete, issues which play an important role in athletic performance.

We have also provided different models which include several delays or an integro-differential approach.

% \appendix

% \section{}
% \label{}

\end{document}